\documentclass{article}

\usepackage[colorlinks]{hyperref}
\usepackage{amssymb}
\usepackage[centertags]{amsmath}
\usepackage[margin=1in]{geometry}
\usepackage{amsthm}
\usepackage{breqn}
\usepackage{setspace}
\usepackage{graphicx}
\usepackage{enumitem}
\usepackage{float}
\usepackage{subfigure}
\usepackage{amsmath, amssymb, graphics, setspace}

\newcommand{\mathsym}[1]{{}}
\newcommand{\unicode}[1]{{}}

\begin{document}

\title{Bifurcation Analysis of Predator-Prey System using Conformable Fractional Order Discretization}
\maketitle
\begin{center}
\author{${}^{}$ Muhammad Rafaqat$^{1}$,  Abubakar Masha$^{2}$, Nauman Ahmed$^{1,3}$,  Ali Raza$^{4}$, Wojciech Sumelka$^{5}$}
\end{center}
\begin{center}
${}^{1}$ Department of Mathematics and Statistics, University of Lahore, Pakistan\\ 
${}^{2}$ Department of Mathematics and Computer Science, Borno State University, Nigeria \\ 
${}^{3}$ Department of Computer Science and Mathematics, Lebanese American University, Beirut, Lebanon
${}^{4}$ Department of Physical Analysis, University of Chenab, Gujrat, Pakistan
${}^{5}$ Institute of Structural Analysis, Poznan University of Technology, Poznan, Poland
\end{center}

\begin{center}
E.mail: 
${}^{5}$ wojciech.sumelka@put.poznan.pl
\end{center}

\maketitle

\numberwithin{equation}{section}
\newtheorem{theorem}{Theorem}[section]
\newtheorem{proposition}[theorem]{Proposition}
\newtheorem{remark}[theorem]{Remark}
\newtheorem{corollary}[theorem]{Corollary}
\newtheorem{lemma}[theorem]{Lemma}
\newtheorem{definition}[theorem]{Definition}
\newtheorem{example}[theorem]{Example}
\newtheorem{conjecture}{Conjecture}

\section{Abstract}

In this paper, conformal fractional order discretization \cite{35,s1,s2} is used to analyze bifurcation analysis and stability of a predator-prey system. A continuous model has been discretized into a discrete one while preserving the fractional-order dynamics. This allows us to look more closely at the stability properties of the system and bifurcation phenomena, including period-doubling and Neimark-Sacker bifurcation. Through numerical and theoretical methods, this research investigated how the modification in system parameters affects the overall dynamics, which may have implications for ecological management and conservation strategies.
\newline
\textbf{Keywords.} Predator-Prey model; Fractional differential equations; Bifurcation; Stability results; simulations\\
\section{Introduction}

The fractional-order predator-prey system \cite{1,2} is a development of the classical predator-prey model, such as the Lotka-Voltera system that describes the interaction between two species of predator and prey. In conventional models, ordinary differential equations with derivatives of integer order are typically used to illustrate population dynamics. The use of fractional calculus that deals with the derivatives and integrals of non-integer order gives rise to the fractional order predator-prey system \cite{4}. The classical Lotka-Volterra predator-prey model is represented by the following set of ordinary differential equations (ODEs):\\
\begin{equation}\label{R1}
\begin{aligned}
\frac{dP}{dt} &= \alpha P - \beta P H, \\
\frac{dH}{dt} &= \delta P H - \gamma H,
\end{aligned}
\end{equation}

where:
\begin{itemize}
    \item \(P\) represents the population of prey (e.g., rabbits),
    \item \(H\) represents the population of predators (e.g., foxes),
    \item \(\alpha\), \(\beta\), \(\delta\), and \(\gamma\) are parameters representing the birth rate of prey, predation rate, natural death rate of prey, and death rate of predators, respectively.
\end{itemize}

To introduce fractional-order derivatives into this model, one replaces the integer-order derivatives with fractional-order derivatives. For example, the fractional-order predator-prey model can be represented as follows:

\begin{equation}\label{R2}
\begin{aligned}
\frac{d^\alpha P}{dt^\alpha} &= \alpha P - \beta P H, \\
\frac{d^\alpha H}{dt^\alpha} &= \delta P H - \gamma H,
\end{aligned}
\end{equation}

Here, \( \frac{d^\alpha}{dt^\alpha} \) denotes the fractional derivative of order \( \alpha \), which introduces memory effects and non-local interactions into the system, allowing for a more accurate representation of real-world ecological dynamics.
It has been noted that the specific form of the fractional derivative depends on the chosen fractional calculus approach (e.g., Riemann-Liouville, Caputo, etc.), and the choice of order \( \alpha \) affects the dynamics of the system.\\
Recently, some researchers have explored the applications of fractional order differential equations in various biological scenarios, such as ecological systems with delays \cite{15,16,17}, epidemiological systems rooted in control \cite{18}, and ecological systems involving diffusion \cite{19, 25}, among others. These equations also find utility in a range of scientific and engineering disciplines \cite{20,23,24}. Furthermore, in \cite{21,22} author discusses the Hopf bifurcation and chaos control in a fractional-order modified hybrid optimal system and fractional-order hyperchaotic system. Recent research delves into topics like approximating solutions for non-linear fractional-order differential population models \cite{27,28}. We refer the reader \cite{26} while others investigate the qualitative dynamics of nonlinear interactions within biological systems \cite{29,30}. For more details on the bifurcation we refer the reader to \cite{31,32}.\\
This study investigates the stability analysis of a discrete-time predator-prey model using fractional-order discretization, clarifying the complexities of ecological dynamics through a lens that transcends traditional boundaries. By utilizing the power of fractional calculus, the study aims to shed more light on the intricate dynamics of predator-prey interactions, revealing new perspectives and enriching the understanding of ecological systems. Through rigorous analysis and simulation, the study seeks to clarify the underlying mechanisms governing system stability and resilience, paving the way for more informed conservation and management strategies in the face of environmental change.\\

In this paper, we study the following Lotka-Volterra systems
\begin{equation}\label{R3}
\begin{cases}
x_{n+1}= x_n e^{\left(r (1 -x_n)-b y_n\right)\dfrac{h^\alpha}{\alpha}}\\
y_{n+1}= y_n e^{\left(b x_n -d\right)\dfrac{h^\alpha}{\alpha}}\\
\end{cases}
\end{equation}
where $r,b,d,h>0$ and $0<\alpha\leq1$
The Jacobian of the system \eqref{R3} is\\
$J(E^*)=$
\begin{equation}\label{RJ1}
\left(
\begin{array}{cc}
 e^{-\dfrac{h^{\alpha } (r (-1+x_n)+b y_n)}{\alpha }}\left(\dfrac{-h^{\alpha } r x_n +\alpha}{\alpha}\right) & -e^{\dfrac{-h^{\alpha}(r (-1+x_n)+by_n)}{\alpha }} \left(\dfrac{bh^{\alpha } x_n}{\alpha}\right)\\
 e^{\dfrac{h^{\alpha } (-d+b x_n)}{\alpha }} \left(\dfrac{b h^{\alpha } y_n}{\alpha }\right) & e^{\dfrac{h^{\alpha } (-d+b x_n)}{\alpha }} \\
\end{array}
\right)
\end{equation}

\ Section {Existance and Stability of fixed points}
\begin{lemma}\label{RL1}
The system \eqref{R3} has the following positive fixed points if $0 <\alpha\leq1$ and $b>d$
$$E^*=\left(\dfrac{d}{b}, \dfrac{(b-d) r}{b^2}\right)$$
\end{lemma}
The Jacobian of the system \eqref{R3} at any point $\left(\dfrac{d}{b}, \dfrac{(b -d) r}{b^2}\right)$ can be evaluated as
\begin{equation}\label{RJ2}
\left(
\begin{array}{cc}
1-\dfrac{d r h^{\alpha }}{b \alpha } & -\dfrac{d h^{\alpha }}{\alpha } \\
\dfrac{h^{\alpha } (b-d) r}{b \alpha } & 1 \\
\end{array}
\right)\\
\end{equation}

Next, to discuss the stability of the fixed point of $\eqref{R3}$, we follow Lemma 1 from \cite{s3,24},
therefore the characteristic equation of the Jacobian matrix \eqref{RJ2} is the following:
\begin{equation}\label{CH1}
\lambda^2 +p \lambda +q =0
\end{equation}
$p =-2+\dfrac{d r h^{\alpha }}{b \alpha}$\\
$q = \dfrac{-d r h^{\alpha } (d h^{\alpha } +\alpha )+b \left({h^{2\alpha }} d r +\alpha ^2\right)}{b \alpha ^2}$

\begin{theorem}\label{th1}
The following conditions hold about the stability of the positive fixed point $E^*$
\begin{enumerate}
\item The positive fixed point $E^*$ is sink if $d<b$ and $b h^{\alpha }<h^{\alpha } d+\alpha$ and one of the following condition holds
\begin{itemize}
\item $h^{\alpha }\geq \dfrac{2 \alpha }{b-d}$

\item $r <\dfrac{4 b \alpha ^2}{c d (-b h^{\alpha }+h^{\alpha } d+2 \alpha )}$
\end{itemize}

\item The positive fixed point $E^*$ is Saddle if the following condition holds
\begin{itemize}
\item $d<b$, $b h^{\alpha }<h^{\alpha } d+2 \alpha$ and $ b \left(h^{2\alpha } d r+4 \alpha ^2\right)<h^{\alpha } d r (h^{\alpha } d +2 \alpha )$
\end{itemize}

\item The positive fixed point $E^*$ is Source if $d <b$, $h^{\alpha } d +\alpha < b h^{\alpha }$ and one of the following condition are hold
\begin{itemize}
\item $h^{\alpha }\geq \dfrac{2 \alpha }{b-d}$

\item $r <\dfrac{4 b \alpha ^2}{h^{\alpha } d (-b h^{\alpha } +d h^{\alpha }+2 \alpha )}$
\end{itemize}

\item E* is non-hyperbolic and observed a flip bifurcation if the following conditions are hold
\begin{itemize}
\item $h^{\alpha }>\dfrac{2 \alpha }{r}$ and $b=\dfrac{h^{\alpha } d r (h^{\alpha } d +2 \alpha )}{d h^{2\alpha}r +4 \alpha ^2}$

\item $\dfrac{h^{\alpha } d r}{b \alpha }\neq 2,4$
\end{itemize}

\item The positive fixed points $E^*$ are complex if the following conditions are hold
\begin{itemize}
\item $0 <\dfrac{h^{\alpha } d r}{b \alpha} < 4$

\item $b = d + \dfrac{\alpha}{h^{\alpha }}$
\end{itemize}

\end{enumerate}
\end{theorem}

\section{Neimark Sacker Bifurcation and Flip Bifurcation}
Now, we investigate the Neimark-Sacker and Flip bifurcation at the positive fixed point $E^*$\\
Let us define a set\\
$A_{NS}= \left\{ (\alpha, b, h, d, r): \alpha, b, h, d, r >0, 0<\alpha \leq 1\text{ and }b =d +\dfrac{\alpha }{h^{\alpha }}\text{ and } 0 <\dfrac{h^{\alpha } d r}{b \alpha} < 4\right\}$\\
Let $b$ be a bifurcation parameter. The Neimark-Sacker bifurcations vary in a small neighborhood of $A_{NS}$.\\
Choose $\epsilon$ as a small perturbation parameter, where $(|\epsilon|\ll 1)$, then the perturbation of system \eqref{R3} is as follows:

\begin{equation}\label{R4}
\begin{cases}
x_{n+1}= x_n e^{\left(r (1 -x_n)-(b+\epsilon) y_n\right)\dfrac{h^\alpha}{\alpha}}\\
y_{n+1}= y_n e^{(b x_{n} -d)\dfrac{h^\alpha}{\alpha}}\\
\end{cases}
\end{equation}
where $r,b,d,h>0$ and $0<\alpha\leq1$\\
Let us define a translation $u_n = x_n-\dfrac{d}{b}$ and  $v_n= y_n-\dfrac{(b-d) r}{b^2}$ that translate the positive equilibrium point of the system \eqref{R3} to the origin and we have the following
\begin{equation}\label{R5}
\begin{cases}
u_{n+1}=C_1u_n + C_2v_n +C_3 u^2_n +C_4v^2_n +C_5u_nv_n +C_6u^2_nv_n +C_7 u_nv^2_n + C_8u^3_n + C_9v^3_n +O((|u_n|+|v_n|)^4)\\
v_{n+1}= D_1u_n +D_2v_n +D_3u^2_n+D_5u_nv_n + D_6u^2_nv_n +D_8u^3_n+O((|u_n|+|v_n|)^4)\\
\end{cases}
\end{equation}
where ${C_{i}'s}$ and ${D_{i}'s}$ are given below\\
$C_1 = \dfrac{-h^{2\alpha } d r+ d \alpha h^{\alpha } +\alpha ^2+ \alpha h^{\alpha } \epsilon }{\alpha ( d h^{\alpha }+\alpha + \epsilon h^{\alpha } )}$,\,
$C_2 = -\dfrac{h^{\alpha } d}{\alpha}$,\,
$C_3 = \dfrac{r h^{\alpha } \left(h^{2\alpha } d r-2 d \alpha h^{\alpha } -2 \alpha ^2-2  \alpha  \epsilon h^{\alpha } \right)}{2 \alpha ^2 ( d h^{\alpha }+\alpha +\epsilon h^{\alpha } )}$,\\
$C_4 = \dfrac{ d h^{\alpha } ( d h^{\alpha }+\alpha + \epsilon h^{\alpha } )}{2 \alpha ^2}$,\,
$C_5 = \dfrac{-h^{2\alpha} d r+ d \alpha h^{\alpha }+\alpha ^2+ \alpha  \epsilon h^{\alpha } }{\alpha ^2}$,\\
$C_6 = -\dfrac{ r h^{\alpha } \left(h^{2\alpha } d r-2 d \alpha h^{\alpha } -2 \alpha ^2-2 \alpha  \epsilon h^{\alpha }\right)}{2 \alpha^3}$,
$C_7 = \dfrac{( d h^{\alpha }+\alpha + \epsilon h^{\alpha }) \left(-h^{2\alpha } d r+ d \alpha h^{\alpha } +\alpha ^2+ \alpha  \epsilon h^{\alpha }\right)}{2 \alpha^3},$\\
$C_8 = -\dfrac{\left(h^{2\alpha } r^2 \left(h^{2\alpha } d r-3  d \alpha h^{\alpha }-3 \alpha ^2-3  \alpha  \epsilon h^{\alpha } \right)\right)}{6\alpha^3 ( d h^{\alpha }+\alpha + \epsilon h^{\alpha })}$,
$C_9 = -\dfrac{\left(d h^{\alpha } ( d h^{\alpha }+\alpha + h^{\alpha } \epsilon  )^2\right)}{6 \alpha^3}$\\
$D_1 = \dfrac{h^{\alpha } r (\alpha +h^{\alpha } \epsilon )}{\alpha  (h^{\alpha } d+\alpha +h^{\alpha } \epsilon )}$,
$D_2 = 1$,
$D_3 = \dfrac{h^{\alpha } r (\alpha +h^{\alpha } \epsilon )}{2 \alpha ^2}$,
$D_5 = \dfrac{(h^{\alpha } d+\alpha +h^{\alpha } \epsilon )}{\alpha}$\\
$D_6 = \dfrac{(h^{\alpha } d +\alpha +h^{\alpha } \epsilon )^2}{2 \alpha ^2}$,
$D_8 = \dfrac{h^{\alpha } r (\alpha +h^{\alpha } \epsilon ) (h^{\alpha } d+\alpha +h^{\alpha } \epsilon )}{6 \alpha ^3}$,\\
Let,\\
\begin{equation}\label{R6}
\lambda^2- p(\epsilon)\lambda+q(\lambda)=0
\end{equation}
be the characteristic equation of the above system, where $p(\epsilon)$ and $q(\epsilon)$ are given by the following\\
$p(\epsilon)=-\dfrac{h^{2\alpha } d r-2 \alpha ^2-2 ch^{\alpha } \alpha  (d+\epsilon )}{\alpha  (\alpha +h^{\alpha } (d+\epsilon ))}$\\
$q(\epsilon)=\dfrac{\alpha ^3+h^{3\alpha }d r \epsilon +h^{\alpha } \alpha ^2 (d+\epsilon )}{\alpha ^2 (\alpha +h^{\alpha } (d+\epsilon ))}$\\
since $(\alpha, b, c, d, r,h)\in A_{NS}$, then $|\lambda_1|=1$ and $|\lambda_2|=1$, where $\lambda_1$ and $\lambda_2$ are solutions
of \eqref{R6} and are given by
$$\lambda_{1,2}=\dfrac{p(\epsilon)}{2}\pm i \dfrac{\sqrt{4q(\epsilon)-p(\epsilon)^2}}{2}$$
Since $|\lambda_1|=|\lambda_2|=\sqrt{q(\epsilon)}$ and
at $\epsilon=0$, we have the following result\\
$\dfrac{d|\lambda_1|}{d\epsilon}=\dfrac{d|\lambda_2|}{d\epsilon}=\dfrac{h^{3\alpha} d r}{2 \alpha ^2 (h^{\alpha } d+\alpha )}\neq 0$ and\\
$p(0)=\dfrac{-h^{2\alpha } d r+2 h^{\alpha } d \alpha +2 \alpha ^2}{\alpha  (h^{\alpha } d+\alpha )}\neq -2,2$\\
 and moreover
\begin{enumerate}
\item $p(0)\neq 1$ if   $d \neq \dfrac{\alpha^2}{h^{2\alpha } r - h^{\alpha } \alpha}$
\item $p(0)\neq 0$ if $d \neq \dfrac{2\alpha^2}{h^{2\alpha } r - 2h^{\alpha } \alpha}$
\end{enumerate}
which is equivalent to say that $\lambda^m, \overline{\lambda}^m\neq1, (m=1,2,3,4)$ if the above two conditions hold.\\
Let us define a transformation \\
$\left(
  \begin{array}{c}
    u_n \\
    v_n \\
  \end{array}
\right)
= \left(
\begin{array}{cc}
 -\dfrac{h^{\alpha } d}{\alpha } & 0 \\
 \dfrac{h^{2\alpha } d r}{2 h^{\alpha } d \alpha +2 \alpha ^2} & -\dfrac{1}{2} \sqrt{\dfrac{h^{2\alpha} d r \left(-h^{2\alpha } d r+4 h^{\alpha } d \alpha +4 \alpha ^2\right)}{\alpha ^2 (h^{\alpha } d+\alpha
)^2}} \\
\end{array}
\right)
\left(
  \begin{array}{c}
    e_n \\
    f_n \\
  \end{array}
\right)$\\
 we get the following canonical form of the system \eqref{R5}\\
 $\left(
  \begin{array}{c}
    e_{n+1} \\
    f_{n+1} \\
  \end{array}
\right)
=\left(
\begin{array}{cc}
 \dfrac{-h^{2\alpha } d r+2 h^{\alpha } d \alpha +2 \alpha ^2}{2 \alpha  (h^{\alpha } d+\alpha )} & -\dfrac{1}{2} \sqrt{\dfrac{(h^{\alpha })^2 d r \left(-(h^{\alpha })^2 d r+4 h^{\alpha } d \alpha +4 \alpha ^2\right)}{\alpha
^2 (h^{\alpha } d+\alpha )^2}} \\
 \dfrac{1}{2} \sqrt{\dfrac{h^{2\alpha } d r \left(-(h^{\alpha })^2 d r+4 h^{\alpha } d \alpha +4 \alpha ^2\right)}{\alpha ^2 (h^{\alpha } d+\alpha )^2}} & \dfrac{-(h^{\alpha })^2 d r+2 h^{\alpha } d \alpha +2 \alpha
^2}{2 \alpha  (h^{\alpha } d+\alpha )} \\
\end{array}
\right)\\
\left(
  \begin{array}{c}
    e_n \\
    f_n \\
  \end{array}
\right)
+
 \left(
  \begin{array}{c}
    F(e_n, f_n) \\
    G(e_n, f_n) \\
  \end{array}
\right)$\\
Where $F(e_n, f_n)$ and $G(e_n, f_n)$ are given by
\begin{equation}\label{R7}
\begin{cases}
F(e_n, f_n)= C_{11} e^2_n+C_{12} f^2_n+ C_{13} e_n f_n+ C_{14} e^2_n f_n+ C_{15} e_n f^2_n +C_{16} e^3_n+ C_{17}f^3_n+O((|e_n|+|f_n|)^4)\\
G(e_n, f_n)=D_{11} e^2_n+D_{12} f^2_n+ D_{13} e_n f_n+ D_{14} e^2_n f_n+ D_{15} e_n f^2_n+D_{16} e^3_n+ D_{17}f^3_n +O((|e_n|+|f_n|)^4)\\
\end{cases}
\end{equation}
and $C_{ij}'s$ and $D_{i,j}'s$ are given by\\
\begin{eqnarray*}
C_{11}&=&-\dfrac{h^{2\alpha } d r \left(h^{2\alpha } d r-4 h^{\alpha } d \alpha -4 \alpha ^2\right)}{8 \left(\alpha ^3 (h^{\alpha } d+\alpha )\right)},\\
C_{12}&=& \dfrac{h^{2\alpha } d r \left(h^{2\alpha } d r-4 h^{\alpha } d \alpha -4 \alpha ^2\right)}{8 \alpha ^3 (h^{\alpha } d+\alpha )},\\
C_{13}&=&  -\dfrac{h^{2\alpha } d r-2 h^{\alpha } d \alpha -2 \alpha ^2}{2 (\alpha  (h^{\alpha } d+\alpha ))}-\dfrac{\left(\sqrt{-\dfrac{h^{2\alpha } d r \left(h^{2\alpha } d r-4 h^{\alpha } d \alpha -4 \alpha ^2\right)}{\alpha ^2 (h^{\alpha } d+\alpha )^2}} \left(h^{2\alpha } d r-2 h^{\alpha } d \alpha -2 \alpha ^2\right)\right) }{4 \alpha ^2},\\
C_{14}&=& -\dfrac{h^{2\alpha } d r \left(h^{2\alpha } d r-4 h^{\alpha } d \alpha -4 \alpha ^2\right)}{8\alpha ^3 (h^{\alpha } d+\alpha )}-\dfrac{\left(h^{2\alpha } d r \left(h^{2\alpha } d r-4 h^{\alpha } d \alpha -4 \alpha^2\right) \sqrt{-\dfrac{h^{2\alpha } d r \left(h^{2\alpha } d r-4 h^{\alpha } d \alpha -4 \alpha ^2\right)}{\alpha ^2 (h^{\alpha } d+\alpha )^2}}\right) }{16\alpha^4},\\
C_{15}&=& \dfrac{h^{2\alpha } d r \left(h^{2\alpha } d r-4 h^{\alpha } d \alpha -4 \alpha ^2\right) \left(h^{2\alpha } d r-2 h^{\alpha } d \alpha -2 \alpha ^2\right) }{16 \alpha ^5 (h^{\alpha } d+\alpha
)},\\
C_{16}&=& -\dfrac{h^{4\alpha } d^2 r^2 \left(h^{2\alpha } d r-6 h^{\alpha } d \alpha -6 \alpha ^2\right)}{48 \left(\alpha ^5 (h^{\alpha } d+\alpha )\right)},\\
C_{17}&=& \dfrac{h^{2\alpha } d r \left(h^{2\alpha } d r-4 h^{\alpha } d \alpha -4 \alpha^2\right) \sqrt{-\dfrac{h^{2\alpha } d r \left(h^{2\alpha } d r-4 h^{\alpha } d \alpha -4 \alpha ^2\right)}{\alpha^2 (h^{\alpha } d+\alpha )^2}} }{48 \alpha ^4},\\
D_{11}&=& -\dfrac{h^{4\alpha } d^2 r^2 \left(h^{2\alpha } d r-4 h^{\alpha } d \alpha -4 \alpha ^2\right)}{8 \left(\alpha ^4 (h^{\alpha } d+\alpha )^2 \sqrt{-\dfrac{h^{2\alpha } d r \left(h^{2\alpha } d r-4 h^{\alpha } d \alpha -4 \alpha ^2\right)}{\alpha ^2 (h^{\alpha } d+\alpha )^2}}\right)},\\
D_{12}&=& \dfrac{h^{4\alpha } d^2 r^2 \left(h^{2\alpha } d r-4 h^{\alpha } d \alpha -4 \alpha ^2\right) }{8 \alpha ^4 (h^{\alpha } d+\alpha )^2 \sqrt{\dfrac{h^{2\alpha } d r \left(-h^{2\alpha } dr +4 h^{\alpha } d \alpha +4 \alpha ^2\right)}{\alpha ^2 (h^{\alpha } d+\alpha )^2}}},\\
D_{13}&=& \left(-\dfrac{h^{\alpha } d (h^{\alpha } d+\alpha )}{\alpha ^2}-\dfrac{h^{2\alpha } d r \left(h^{2\alpha } d r-2 h^{\alpha } d \alpha -2 \alpha ^2\right)}{4 \alpha ^3 (h^{\alpha } d+\alpha)}\right),\\
D_{14}&=& \left(-\dfrac{h^{4\alpha } d^2 r^2 \left(h^{2\alpha } d r-4 h^{\alpha } d \alpha -4 \alpha ^2\right)}{16 \alpha ^5 (h^{\alpha } d+\alpha )}+\dfrac{d^2 \left(h^{4\alpha } d^2+2 h^{3\alpha }d \alpha +h^{2\alpha } \alpha ^2\right)}{2 \alpha ^4}\right),\\
\end{eqnarray*}
\begin{eqnarray*}
D_{15}&=& \dfrac{h^{4\alpha } d^2 r^2 \left(h^{2\alpha } d r-4 h^{\alpha } d \alpha -4 \alpha ^2\right) \left(h^{2\alpha } d r-2 h^{\alpha } d \alpha -2 \alpha ^2\right) }{16 \alpha ^6 (h^{\alpha }d+\alpha )^2 \sqrt{-\dfrac{h^{2\alpha } d r \left(h^{2\alpha } d r-4 h^{\alpha } d \alpha -4 \alpha ^2\right)}{\alpha ^2 (h^{\alpha } d +\alpha )^2}}},\\
D_{16}&=&  -\dfrac{d^3 \left(h^{8\alpha } d r^4+8h^{7\alpha } d^3 r \alpha -6 h^{7\alpha } d r^3 \alpha +24 h^{6\alpha } d^2 r \alpha ^2-6 h^{6\alpha } r^3 \alpha ^2 +24 h^{5\alpha } d r \alpha^3+8 h^{4\alpha } r \alpha ^4\right)}{48 \left(\alpha ^6 (h^{\alpha } d+\alpha )^2 \sqrt{-\dfrac{h^{2\alpha } d r \left(h^{2\alpha } d r-4 h^{\alpha } d \alpha -4 \alpha ^2\right)}{\alpha ^2 (h^{\alpha } d+\alpha )^2}}\right)},\\
D_{17}&=&  \dfrac{h^{4\alpha } d^2 r^2 \left(h^{2\alpha } d r-4 h^{\alpha } d \alpha -4 \alpha ^2\right) }{48 \alpha ^5 (h^{\alpha } d+\alpha )},
\end{eqnarray*}
To examine the Neimark-Sacker bifurcation, let us consider the first Lyapunov exponent, derived as follows:
$$a= \left[-Re\left(\frac{(1-2\lambda)\overline{\lambda}^2}{1-\lambda} \xi_{20}\xi_{11}-\frac{1}{2}|\xi_{11}|^2-|\xi_{02}|^2+Re(\overline{\lambda}\xi_{21})\right)\right]_{\epsilon=0}$$
where\\
$\xi_{20}=\dfrac{1}{8}\left[(F_{e_ne_n}-F_{f_nf_n}+2G_{e_nf_n})+i(G_{e_ne_n}-G_{f_nf_n}-2F_{e_nf_n})\right].$\\
$\xi_{11}=\dfrac{1}{4}\left[(F_{e_ne_n}+F_{f_nf_n})+i(G_{e_ne_n}+G_{f_nf_n})\right].$\\
$\xi_{02}=\dfrac{1}{8}\left[(F_{e_ne_n}-F_{f_nf_n}-2G_{e_nf_n})+i(G_{e_ne_n}-G_{f_nf_n}+2F_{e_nf_n})\right].$\\
$\xi_{21}=\dfrac{1}{16}\left[(F_{e_ne_ne_n}+F_{e_nf_nf_n}+G_{e_ne_nf_n}+G_{f_nf_nf_n})+i(G_{e_ne_ne_n}+G_{e_nf_nf_n}-F_{e_ne_nf_n}-F_{f_nf_nf_n})\right].$
\begin{theorem}
Let $(\alpha, b, h, d, r)\in A_{NS}$ and $b\ne0$, then the system \eqref{R3} observes a Neimark-Sacker bifurcation at the positive fixed point $E^{*}$ when the bifurcation parameter varies in a small neighbourhood of $b$. Furthermore, when $b$ is negative, an attracting invariant curve emanates from $E^*$. When $b$ is positive, a repelling invariant curve emanates from $E^{*}$.
\end{theorem}
Let us define a set\\
$A_{PD}=\left\{(\alpha, b,d,h): \alpha, b,d,h>0, 0<\alpha\le 1, h^{\alpha }>\dfrac{2 \alpha }{r}, b=\dfrac{h^{\alpha } d r (h^{\alpha } d +2 \alpha )}{d {h^{2\alpha }}  r +4 \alpha ^2}\,\,
\text{and}\,\, \dfrac{h^{\alpha } d r}{b \alpha }\neq 4\right\}$\\
The system \eqref{R3} observes a period-doubling bifurcation within a domain $A_{PD}$. Let $b$ be a bifurcation parameter and $\epsilon^* (|\epsilon^*|\ll 1)$ be a perturbation parameter, then the perturbation of the system \eqref{R3} is the following:

\begin{equation}\label{R8}
\begin{cases}
x_{n+1}= x_n e^{(r (1 -x_n)-(b+\epsilon^*) y_n)\dfrac{h^\alpha}{\alpha}}\\
y_{n+1}= y_n e^{((b +\epsilon^*) x_n -d)\dfrac{h^\alpha}{\alpha}}\\
\end{cases}
\end{equation}

Let $u_n=x_n- \dfrac{d}{b}$ and $v_n=y_n-\dfrac{(b -d) r}{b^2}$, then the positive fixed point $E^{*}$ transform the system $\eqref{R8}$ into origin and we get the following:\\
\begin{equation}\label{R9}
\left(
  \begin{array}{c}
    u_{n+1} \\
    v_{n+1} \\
  \end{array}
\right)
=
\left(
  \begin{array}{c}
    c_{11}u_n+ c_{12}v_n+c_{13}u^2_n+c_{14}v^2_n+c_{15}u_nv_n+c_{16}u_n\epsilon^*+c_{18}u^3_n+\\
    c_{20}u^2_nv_n+ c_{21}u^2_n\epsilon^*+ c_{22}u_nv_n\epsilon^*+O((|u_n|+|v_n|+\epsilon^*)^4)\\
    d_{11}u_n+ d_{12}v_n +d_{13}u^2_n +d_{15}u_nv_n +d_{18}u^3_n +d_{20}u^2_nv_n+\\
    d_{21}u^2_n\epsilon^*+d_{22}u_nv_n\epsilon^*+O((|u_n|+|v_n|+\epsilon^*)^4)
  \end{array}
\right)
\end{equation}

where $c_{ij} \text{and} d_{ij}$ are given by\\
\begin{eqnarray*}
c_{11}&=&\dfrac{-h^{2\alpha } d r +h^{\alpha } d \alpha -2 \alpha ^2}{\alpha  ( dh^{\alpha } +2 \alpha )},\,\,
c_{12}=-\dfrac{ d h^{\alpha } }{\alpha }, \,\,
c_{13}= \dfrac{d \left(h^{3\alpha } r^2 -2h^{2\alpha } r \alpha \right)}{2 \alpha ^2 ( d h^{\alpha } +2 \alpha )},\\
c_{14}&=& \dfrac{h^{3\alpha } d^2 ( d r h^{\alpha } +2 r \alpha )}{2 \alpha ^2 \left(h^{2\alpha } d r+4 \alpha ^2\right)},\,\,
c_{15}=\dfrac{h^{2\alpha } \left(h^{2\alpha } d^2 r^2 -(h^{\alpha }) d^2 r \alpha +2 d r \alpha ^2\right)}{\alpha ^2 \left(h^{2\alpha } d r +4 \alpha ^2\right)},\\
c_{16}&=&\dfrac{\left(h^{4\alpha } d^2 r^2+8h^{2\alpha } d r \alpha ^2+16 \alpha ^4\right) }{ d r h^{\alpha }\alpha  (d h^{\alpha } +2 \alpha )^2},\\
c_{18}&=&-\dfrac{h^{2\alpha } r^2 \left(h^{2\alpha } d r-3 d h^{\alpha } \alpha -2 \alpha ^2\right)}{6 \left(\alpha ^3 (d h^{\alpha } +2 \alpha )\right)},\,\,
c_{20}=-\dfrac{h^{4\alpha } r^2 \left(h^{\alpha }d^2 r -2 d^2 \alpha \right)}{2 \left(\alpha ^3 \left(h^{2\alpha } d r+4 \alpha ^2\right)\right)}\\
c_{21}&=&\dfrac{\left(-h^{4\alpha } d^2 r^2-8 h^{2\alpha } d r \alpha ^2 -16 \alpha ^4\right) }{2 d \alpha ^2 (d h^{\alpha } +2 \alpha )^2},\,\,
c_{22}= -\dfrac{h^{\alpha } }{\alpha },\,\,
d_{11}=\dfrac{2 (h^{\alpha } r-2 \alpha )}{h^{\alpha } d +2 \alpha }\\
d_{12}&=&1,\,\,
d_{13}=\dfrac{h^{2\alpha } d r (c r -2 \alpha )}{\alpha  \left(h^{2\alpha } d r +4 \alpha ^2\right)},\,\,
d_{15}=\dfrac{h^{2\alpha } d r (h^{\alpha } d +2 \alpha )}{\alpha  \left(h^{2\alpha } d r +4 \alpha ^2\right)},\\
d_{18}&=&\dfrac{h^{4\alpha } d^2 r^2 (h^{\alpha } r-2 \alpha ) (h^{\alpha } d+2 \alpha )}{3 \alpha ^2 \left(h^{2\alpha } d r +4 \alpha ^2\right)^2},\,\,
d_{20}=\dfrac{h^{4\alpha } d^2 r^2 (h^{\alpha } d +2 \alpha )^2}{2 \alpha ^2 \left(h^{2\alpha } d r +4 \alpha ^2\right)^2},\,\,
d_{21}=\dfrac{h^{2\alpha } r }{2 \alpha ^2},\,\,
d_{22}=\dfrac{h^{\alpha }}{\alpha}
\end{eqnarray*}
Construct an invertible matrix
\begin{equation*}
T=
\left(
\begin{array}{cc}
 -\dfrac{h^{\alpha } d +2 \alpha }{h^{\alpha } r -2 \alpha } & -\dfrac{h^{\alpha } d}{2 \alpha } \\
 1 & 1 \\
\end{array}
\right)
\end{equation*}
and define a transformation
\begin{equation*}
\left(
\begin{array}{c}
  e_n \\
  f_n
\end{array}
\right)
=T
\left(
\begin{array}{c}
  u_n \\
  v_n
\end{array}
\right)
\end{equation*}
we get the following\\
\begin{equation}\label{R11}
\left(
\begin{array}{c}
  e_{n+1} \\
  f_{n+1}
\end{array}
\right)=
\left(
\begin{array}{cc}
 -1 & 0 \\
 0 & \dfrac{-h^{2\alpha } d r +3 h^{\alpha } d \alpha +2 \alpha ^2}{\alpha  (h^{\alpha } d +2 \alpha )} \\
\end{array}
\right)
\left(
\begin{array}{c}
  e_{n} \\
  f_{n}
\end{array}
\right)
+
\left(
\begin{array}{c}
  F(e_n, f_n,\epsilon^*) \\
  G(e_n, f_n,\epsilon^*)
\end{array}
\right)
\end{equation}
 where $F(e_n, f_n,\epsilon^*)$  and  $G(e_n, f_n,\epsilon^*)$  are defined by\\
    $F(e_n, f_n,\epsilon^*)=c_1 e^2_n+c_2 e_n f_n+ c_3 f^2_n+c_4\epsilon^*+c_5e_n\epsilon^*+c_6e^2_n\epsilon^*+c_7 f_n\epsilon^*+c_8 e_n f_n\epsilon^*+c_9f^2_n\epsilon^*
    \quad +O((|e_n|+|f_n|+|\epsilon^*|)^4)$\\
    $G(e_n, f_n,\epsilon^*)=d_1 e^2_n +d_2 e_n f_n + d_3 f^2_n+d_4\epsilon^*+d_5e_n\epsilon^*+d_6e^2_n\epsilon^*+d_7 f_n\epsilon^*+d_8 e_n f_n\epsilon^*+b_9f^2_n\epsilon^*\\
    \quad +O((|e_n|+|f_n|+|\epsilon^*|)^4)$\\
where $c_{i}$ and $d_{i}$ are given by\\
\begin{eqnarray*}
c_2&=&\dfrac{h^{4\alpha } d^3 r+4 h^{3\alpha } d^2 r \alpha +4 h^{2\alpha } d r \alpha ^2}{2 \alpha ^2 \left(h^{2\alpha } d r +4 \alpha ^2\right)},\\
c_3&=&\dfrac{h^{5\alpha } d^3 r (d +r) (h^{\alpha } r -2 \alpha )}{4 \alpha ^3 ( d h^{\alpha }+2 \alpha ) \left(h^{2\alpha } d r+4 \alpha ^2\right)}, \,\,\,
c_4=\dfrac{(r h^{\alpha }-2 \alpha ) \left(h^{4\alpha } d^2 r^2+8h^{2\alpha } d r \alpha ^2+16 \alpha ^4\right) }{2 r \alpha ^2 (h^{\alpha } d +2 \alpha ) \left(-h^{2\alpha } d r+4 h^{\alpha }
d \alpha +4 \alpha ^2\right)}\\
c_5&=&\dfrac{\left(h^{4\alpha } d^2 r^2+8 h^{2\alpha } d r \alpha ^2 +16 \alpha ^4\right)}{h^{\alpha } d r \alpha  \left(-h^{2\alpha } d r+4h^{\alpha } d \alpha +4 \alpha ^2\right)}, \,\,\,\,
c_6=\dfrac{(h^{\alpha } d -2 \alpha ) (h^{\alpha } d +2 \alpha )^2 \left(h^{2\alpha } d r +4 \alpha ^2\right) }{2 d (h^{\alpha } r-2 \alpha ) \alpha ^2 \left(h^{2\alpha } d r-4 h^{\alpha } d \alpha
-4 \alpha ^2\right)}\\
c_7&=&\dfrac{(h^{\alpha } r-2 \alpha ) \left(h^{4\alpha } d^2 r^2+8 h^{2\alpha } d r \alpha ^2 +16 \alpha ^4\right) }{2 r \alpha ^2 (h^{\alpha } d+2 \alpha ) \left(-(h^{\alpha })^2 d r+4 h^{\alpha }
d \alpha +4 \alpha ^2\right)}\\
c_8&=&-\dfrac{\left(h^{6\alpha } d^4 r +3 h^{5\alpha } d^3 r \alpha -2 h^{5\alpha } d^2 r^2 \alpha +8 h^{4\alpha } d^2 r \alpha ^2-12 h^{3\alpha } d^2 \alpha ^3-4 h^{3\alpha } d r \alpha ^3 -16h^{2\alpha } d \alpha ^4-16 h^{\alpha } \alpha ^5\right) }{2 \left(\alpha ^3 (h^{\alpha } d+2 \alpha ) \left(-h^{2\alpha } d r+4 h^{\alpha } d \alpha +4 \alpha ^2\right)\right)}\\
c_9&=&(h^{\alpha } r-2 \alpha )(h^{7\alpha } d^5 r+4 h^{6\alpha } d^4 r \alpha -2 h^{6\alpha } d^3 r^2 \alpha -4 h^{5\alpha } d^4 \alpha ^2+4h^{5\alpha } d^3 r \alpha ^2-16 h^{4\alpha } d^2 r \alpha^3\\
&+&48 h^{3\alpha } d^2 \alpha ^4 +32 h^{2\alpha } d \alpha^5)
(8(\alpha ^4 (h^{\alpha } d+2 \alpha )^2(-h^{2\alpha } d r +4 h^{\alpha } d \alpha +4 \alpha ^2)))\\
d_1&=&10, \,\,\,
d_2=-\dfrac{h^{3\alpha } d r (d +r) (h^{\alpha } d +2 \alpha )}{(h^{\alpha } r -2 \alpha ) \alpha  \left((h^{\alpha })^2 d r +4 \alpha ^2\right)},\\
d_3&=&-\dfrac{h^{4\alpha } d^3 r (h^{\alpha } r -2 \alpha )^2}{4 \alpha ^3 (h^{\alpha } d +2 \alpha ) \left(h^{2\alpha } d r +4 \alpha ^2\right)}-\dfrac{h^{3\alpha } d^2 r (h^{\alpha } d +2 \alpha )}{2\alpha ^2 \left(h^{2\alpha } d r +4 \alpha ^2\right)}\\
d_4&=&\dfrac{-h^{2\alpha } d r+3 h^{\alpha } d \alpha +2 \alpha ^2}{\alpha  (h^{\alpha } d +2 \alpha )}-\dfrac{h^{\alpha } \left(h^{4\alpha } d^3 r^2+h^{4\alpha } d^2 r^3+8 h^{2\alpha } d^2 r \alpha ^2+8 h^{2\alpha } d r^2 \alpha ^2 +16 d \alpha ^4 +16 r \alpha^4\right)}{r \alpha  (h^{\alpha } d +2 \alpha )^2 \left(-h^{2\alpha } d r +4 h^{\alpha } d \alpha +4 \alpha ^2\right)},\\
d_5&=&\dfrac{2 \left(h^{4\alpha } d^3 r^2+h^{4\alpha } d^2 r^3+8 h^{2\alpha } d^2 r \alpha ^2+8 h^{2\alpha } d r^2 \alpha ^2+16 d \alpha ^4+16 r \alpha ^4\right) }{d r (h^{\alpha } r-2
\alpha ) (h^{\alpha } d+2 \alpha ) \left(h^{2\alpha } d r-4 h^{\alpha } d \alpha -4 \alpha ^2\right)},\\
d_6&=&-\dfrac{(h^{\alpha } d+2 \alpha ) \left(h^{2\alpha } d r+4 \alpha ^2\right) \left(h^{2\alpha } d^2+2 h^{\alpha } d \alpha -2 h^{\alpha } r \alpha +4 \alpha ^2\right) }{d (h^{\alpha } r-2 \alpha)^2 \alpha  \left(h^{2\alpha } d r-4 h^{\alpha } d \alpha -4 \alpha ^2\right)},\\
d_7&=&-\dfrac{h^{\alpha } \left(h^{4\alpha } d^3 r^2+h^{4\alpha } d^2 r^3+8 h^{2\alpha } d^2 r \alpha ^2+8 h^{2\alpha } d r^2 \alpha ^2+16 d \alpha ^4+16 r \alpha ^4\right) }{r \alpha(h^{\alpha } d+2 \alpha )^2 \left(-h^{2\alpha } d r+4 h^{\alpha } d \alpha +4 \alpha ^2\right)},\\
d_8&=&(-h^{6\alpha } d^4 r+h^{6\alpha } d^2 r^3-5 h^{5\alpha } d^3 r \alpha -5 h^{5\alpha } d^2 r^2 \alpha +2 h^{4\alpha } d^2 r \alpha ^2 +2 h^{4\alpha } d r^2 \alpha ^2-4 h^{3\alpha } d^2\alpha ^3-4 h^{3\alpha } d r \alpha ^3+8 h^{2\alpha } d \alpha ^4\\
&+&8 h^{2\alpha } r \alpha ^4)/((h^{\alpha } r-2 \alpha ) \alpha ^2 (h^{\alpha } d+2 \alpha )(h^{2\alpha } d r-4 h^{\alpha } d \alpha -4 \alpha ^2))\\
d_9&=&-\dfrac{(h^{\alpha } d +2 \alpha )(-h^{4\alpha } d^2 r +4 h^{2\alpha } d \alpha ^2)}{4 \alpha ^3(-h^{2\alpha } d r +4 h^{\alpha } d \alpha +4 \alpha^2)}
-\dfrac{h^{2\alpha } d (h^{\alpha } r -2 \alpha ) \left(h^{4\alpha } d^2 r^2 -8 h^{2\alpha } d^2 \alpha ^2 +8 h^{2\alpha } d r \alpha ^2 -32 h^{\alpha } d \alpha^3 -16 \alpha^4\right)}{4 \alpha ^3 (h^{\alpha } d+2 \alpha )^2 \left(-h^{2\alpha } d r+4 h^{\alpha } d \alpha +4 \alpha ^2\right)}
\end{eqnarray*}
Next we determine the center manifold $W^{c}(0,0,0)$ of the system \eqref{R11} at the fixed point $(0,0)$ in a small neighborhood of $\epsilon^{*}=0$. Let \begin{equation}\label{R12}
W^{c}(0,0)=
\{(e_n, f_n, \epsilon^{*})\in \mathbb{R}:
f_{n}=W(u_{n}, \epsilon^{*})=h_{0}\epsilon^{*}+h_1 e_{n}^2+h_2 e_{n}\epsilon^{*}+h_{3}{\epsilon^{*}}^2 +O(|e_n|+\epsilon^{*})^3
\}
\end{equation}
\begin{eqnarray*}
f_{n+1} &=& W(e_{n+1}, \epsilon^{*})\\
f_{n+1} &=& W(-e_n + F(e_n,f_n, \epsilon^{*}), \epsilon^{*})\\
f_{n+1} &=& W(-e_n +F(e_n, W(e_n,  \epsilon^{*}), \epsilon^{*}), \epsilon^{*})\\
f_{n+1} &=& W(-e_n +F(e_n,h_{0} \epsilon^{*}+ h_1 e^2_n + h_2 e_n \epsilon^{*} +h_3  \epsilon^{*^2} + O(|e_n|+ | \epsilon^{*}|)^3, \epsilon^{*}), \epsilon^{*})\\
f_{n+1} &=& W(-e_n + c_1 e^2_n +c_2 e_n(h_0 \epsilon^{*})) + c_3 h_0^2 \epsilon^{*^2} + c_4 \epsilon^{*} + c_5 e_n  \epsilon^{*} + C_7 h_0  \epsilon^{*^2} O(|e_n| +| \epsilon^{*}|)^3,  \epsilon^{*}\\
f_{n+1} &=& W(-e_n + c_4  \epsilon^{*} + c_1  \epsilon^2_n + (c_5 + c_2 h_0) \epsilon^{*} e_n + h_0(c_3 h_0 + c_7)\epsilon^{*^2} +O(|e_n|,  \epsilon^{*})^3), \epsilon^{*}\\
f_{n+1} &=& h_0  \epsilon^{*} + h_1 (e^2_n + c^2_4 +  \epsilon^{*^2}- 2c_4 e_n  \epsilon^{*}) + h_2 (-e_n \epsilon^{*} + c_4  \epsilon^{*^2}) + h_3 \epsilon^{*^2} + O(|e_n|,  \epsilon^{*})^3 \\
f_{n+1} &=& h_0  \epsilon^{*} +h_1 e^2_n + (h_1 c^2_4 + h_3 + h_2c_4) \epsilon^{*^2} +(2 h_1, c_4 + h_2) e_n  \epsilon^{*} + O(|e_n| +| \epsilon^{*}|)^3
\end{eqnarray*}
Also
\begin{eqnarray*}
f_{n+1} &=& \alpha_1 f_n + G(e_n, f_n,  \epsilon^{*}), \text{where} \,\,
\alpha_1 = \dfrac{-h^{2 \alpha}d q + 3 h^\alpha d \alpha + 2 \alpha^2}{\alpha(h^{\alpha} d + 2 \alpha)}\\
f_{n+1} &=& \alpha_1 f_n + G(e_n, w(e_n \epsilon^{*}),\epsilon^{*}) \\
f_{n+1} &=& \alpha_1 W(e_n,\epsilon^{*}) + G (e_n, h_0 \epsilon^{*} + h_1 e^2_n + h_2 e_n \epsilon^{*} + h_3 \epsilon^{*^2} + O(|e_n|+ \epsilon^{*})^3,\epsilon^{*})\\
f_{n+1} &=& \alpha_1 (h_0 \epsilon^{*} + h_1 e^2_n + h_2 e_n \epsilon^{*} + h_3 \epsilon^{*^2}) + d_1 e^2_n + d_2 e_n(h_0 \epsilon^{*}) + d_3(h^2_0\epsilon^{*^2}) + d_4 \epsilon^{*} + d_5 e_n \epsilon^{*}
 + O(3) +d_7 h_0 \epsilon^{*^2}\\
f_{n+1}&=& (d_4 + \alpha_1 h_2 + d_2 h_0)e_n \epsilon^{*} + (d_7 h_0 + h_3 \alpha_1 + d_3 h^2_0)\epsilon^{*^2} + O(3).
\end{eqnarray*}
compare coefficients we get the following \\
$d_4 + \alpha_1 h_0 = h_0\rightarrow h_0 = \dfrac{d_4}{1 - \alpha_1},$\\
 $h_1 = \alpha_1 h_1 + d_1\rightarrow h_1 = \dfrac{d_1}{1-a_{22}},$\\
$-2h_1c_4 - h_2 = \alpha_1 h_2 +d_2 h_0 + d_5\rightarrow -2h_1c_4 -d_2 h_0 - d_5 = \alpha_1 h_2 +h_2\rightarrow h_2= \dfrac{-2h_1c_4 -d_2 h_0 - d_5}{1+\alpha_1}$\\
$h_3= h_1 c^2_4 + h_2c_4 = d_7h_0 + h_3 \alpha_1 + d_3 h^2_0\rightarrow
h_3 = \dfrac{d_3 h^2_0 + d_7h_0 - h_1c^2_4 - h_2 c_4}{1 - \alpha_1}$\\
Thus the center manifold can be approximated as
\begin{eqnarray*}
f_n = h_0\epsilon^{*} + h_1e^2_n + h_2 e_n \epsilon^{*} + h_3 \epsilon^{*^2} + O (|e_n| +\epsilon^{*})^3
\end{eqnarray*}
where $h_0$, $h_1$, $h_2$ and $h_3$ are given above. Substitute this in \\
$e_{n+1}=- e_n + F(e_n, f_n, \epsilon^{*})$,\\
 we get the following equation that will describe the dynamics of the given system
\begin{eqnarray*}
e_{n+1}&=& -e_n + c_1e^2_n + c_2e_n(h_0\epsilon^{*} + h_1e^2_n + h_2e_n\epsilon^{*} + h_3\epsilon^{*^2}) + c_3 (h_0\epsilon^{*^2} \\
&+& 2h_0\epsilon^{*}e^2_n)+ c_4\epsilon^{*} + c_5 e_n \epsilon^{*} + c_6 e^2_n\epsilon^{*} +c_8e_n\epsilon^{*}(h_0\epsilon^{*}) + c_9\epsilon^{*}(h^2_0 \epsilon^{*^2}) + O(|e_n|+|\epsilon^{*}|)^4\\
e_{n+1}&=& -e_n + c_1 e^2_n + c_4 \epsilon^{*} + c_5 e_n \epsilon^{*} + c_6 e^2_n\epsilon^{*} + c_2 h_0e_n\epsilon^{*} + c_3 h_0e_n\epsilon^{*^2} + c_3 h^2_0\epsilon^{*^2} + c_9 h^2_0\epsilon^{*^3} + c_7h_1e^2_n\epsilon^{*}\\ &+& 2c_3h_0h_1e^2_n\epsilon^{*}
 + c_2h_2e^2_{n}\epsilon^{*} +... O(4)
\end{eqnarray*}
Let us define the right-hand side of the above equation by $"\Phi"$, then the following are holds about the flip bifurcation at (0,0)
\begin{eqnarray*}
\beta_1 &=& \left( \dfrac{\partial^2\Phi}{\partial e_n \partial \epsilon^{*}}+ \dfrac{1}{2} \dfrac{\partial\Phi}{\partial \epsilon^{*}} \dfrac{\partial^2\Phi}{\partial e^2_n} \right) \neq 0\\
\beta_2 &=& \left(\dfrac{1}{6} \dfrac{\partial^3 \Phi}{\partial e^3_n}+ (\dfrac{1}{2} \dfrac{\partial^2\Phi}{\partial e^2_n})^2\right) \neq 0
\end{eqnarray*}
\section{Chaos Control}
In this Section, we apply the hybrid control strategy to control the chaos generated by the period-doubling and Neimark-Sacker bifurcation in the system \eqref{R3}.\\
Consider the following controlled system:
\begin{equation}\label{R13}
\begin{cases}
x_{n+1}= x_n e^{\left(r (1 -x_n)-b y_n\right)\dfrac{h^\alpha}{\alpha}}+(1-\beta)x_n\\
y_{n+1}= y_n e^{\left(b x_n -d\right)\dfrac{h^\alpha}{\alpha}}+(1-\beta)y_n\\
\end{cases}
\end{equation}
where $0<\beta<1$ is the control parameter for the hybrid control method. It is interesting to note that the fixed points for \eqref{R3} and \eqref{R13}
are the same. The Jacobian of the system \eqref{R13} at the positive fixed point $E^{*}$ is given by:
\begin{equation}\label{R14}
J(E^{*})=
\left[
\begin{array}{cc}
 2-\dfrac{d h^{\alpha } r}{b \alpha }-\beta  & -\dfrac{d h^{\alpha }}{\alpha } \\
 \dfrac{(b-d) h^{\alpha } r}{b \alpha } & 2-\beta  \\
\end{array}
\right]
\end{equation}
The trace and determinant of \eqref{R14} are given below:\\
$T=\text{trac}(J)= 4-\dfrac{d h^{\alpha } r}{b \alpha }-2 \beta$\\
$D= \text{det}(J)=\dfrac{(b-d) d h^{2 \alpha } r+\alpha  \left(d h^{\alpha } r+b \alpha  (-2+\beta )\right) (-2+\beta )}{b \alpha ^2}$\\
The Jury condition states that the fixed point $E^{*}$ of the controlled system \eqref{R13} is stable if the following conditions are met:\\
\begin{equation}
\begin{cases}
T+D+1>0\\
-T+D+1>0\\
D-1<0
\end{cases}
\end{equation}
\section{Numerical Examples}
\begin{example}
For fixed values of parameters $\alpha = 0.7, h = 0.3, d = 2, r = 1$ and $b=3.62597$, the system \eqref{R3} has the fixed point $(0.551576, 0.12367)$. The Jacobian matrix of the system \eqref{R3} has the following eigen values $\lambda_1=0.830386 + 0.557189 I$ and $\lambda=0.830386 - 0.557189 I$. Since $|\lambda_1|=1$ and $|\lambda_2|=1$, so the system \eqref{R3} observe a Neimark-Sacker bifurcation. For the bifurcation parameter $b\in[3,4]$ and $x_0=0.5$, $y_0=0.1$, the figure \eqref{B} describes the bifurcation of the system \eqref{R3} and the figure \eqref{BC} describes the bifurcation of the controlled system \eqref{R13}. Figure \eqref{MLE} represents the graph of maximum Lyapunov exponent, and figure \eqref{PP} is the phase portrait for different values of $b$ for the system \eqref{R3}.
\begin{figure}[H]
    \centering
    \subfigure[$(b, x_n)$]{\includegraphics[width=0.45\textwidth]{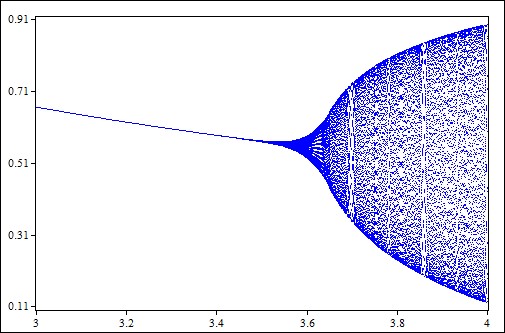}}
    \subfigure[$(b, y_n)$]{\includegraphics[width=0.45\textwidth]{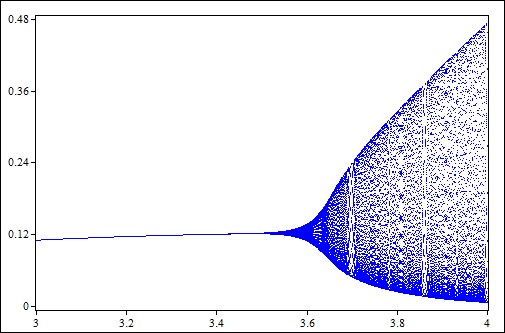}}\\
    \caption{Bifurcation Diagrams for the System\eqref{R3}}
    \label{B}
\end{figure}
\begin{figure}
    \centering
    \subfigure[$(b, x_n), \beta=0.9$]{\includegraphics[width=0.45\textwidth]{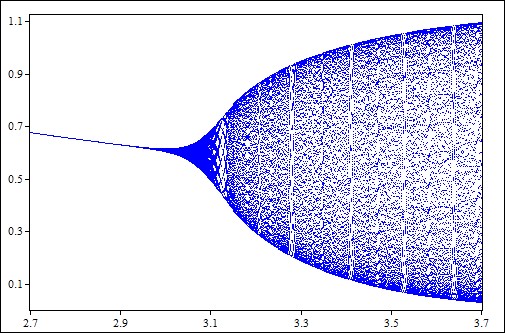}}
    \subfigure[$(b, y_n), \beta=0.9$]{\includegraphics[width=0.45\textwidth]{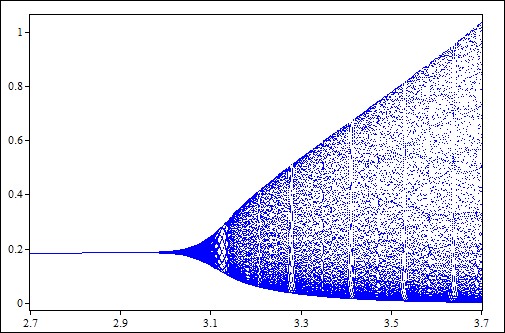}}\\
    \caption{Bifurcation Diagrams for the Controlled System \eqref{R13}}
    \label{BC}
\end{figure}
\newpage
\begin{figure}
    \centering
    \subfigure[$(x_n, y_n)$]{\includegraphics[width=0.50\textwidth]{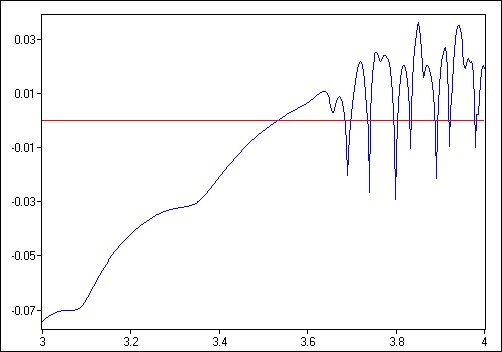}}\\
    \caption{Graph of the Maximum Lyapunov Exponent for the System \eqref{R3}}
    \label{MLE}
\end{figure}
\clearpage
\begin{figure}[H]
    \centering
    \subfigure[$b=3.6$]{\includegraphics[width=0.30\textwidth]{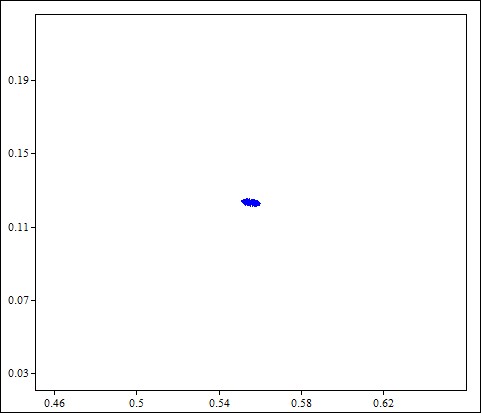}}
    \subfigure[$b=3.61$]{\includegraphics[width=0.30\textwidth]{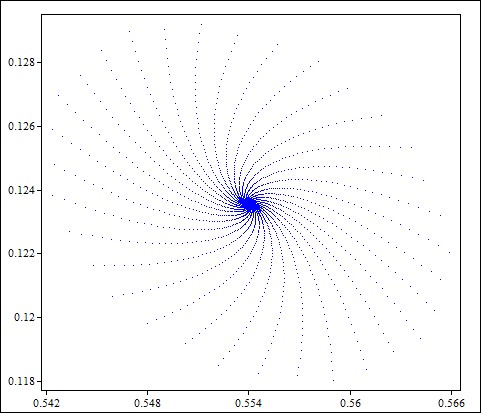}}
    \subfigure[$b=3.62$]{\includegraphics[width=0.30\textwidth]{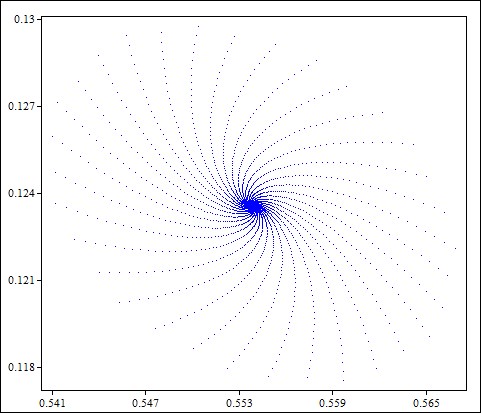}}\\
    \subfigure[$b=3.621$]{\includegraphics[width=0.30\textwidth]{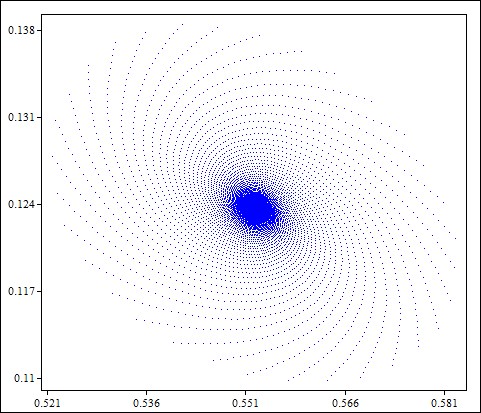}}
    \subfigure[$b=3.622$]{\includegraphics[width=0.30\textwidth]{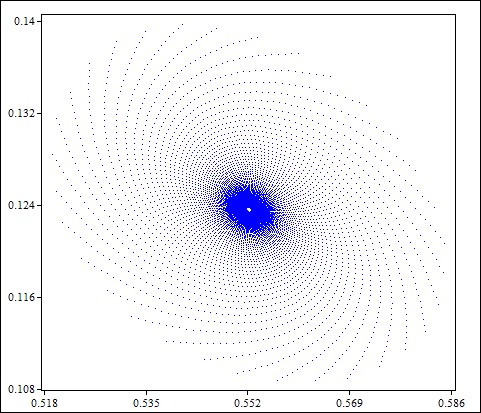}}
    \subfigure[$b=3.623$]{\includegraphics[width=0.30\textwidth]{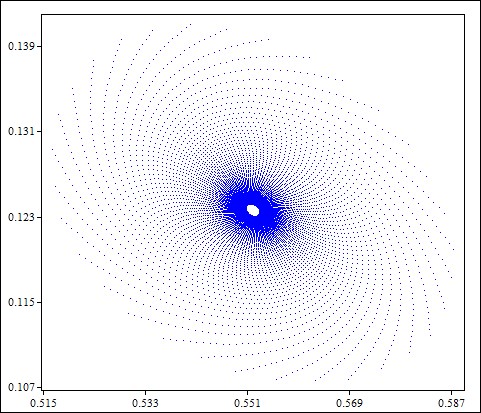}}\\
    \subfigure[$b=3.624$]{\includegraphics[width=0.30\textwidth]{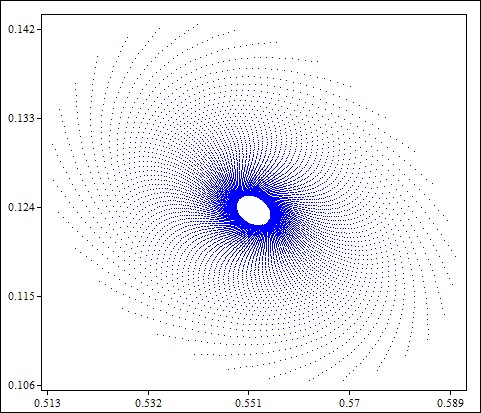}}
    \subfigure[$b=3.625$]{\includegraphics[width=0.30\textwidth]{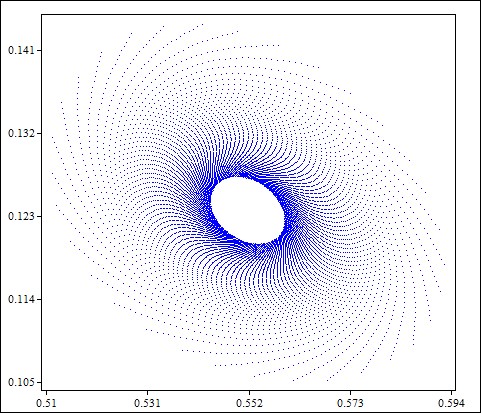}}
    \subfigure[$b=3.626$]{\includegraphics[width=0.30\textwidth]{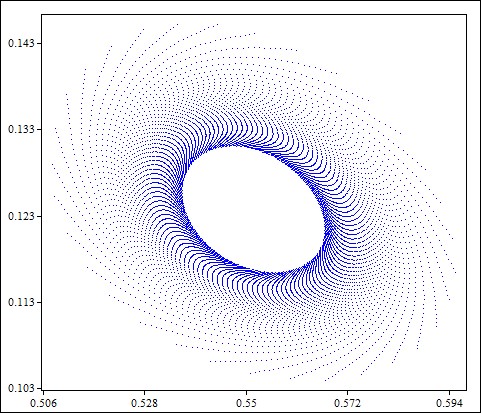}}\\
    \subfigure[$b=3.627$]{\includegraphics[width=0.30\textwidth]{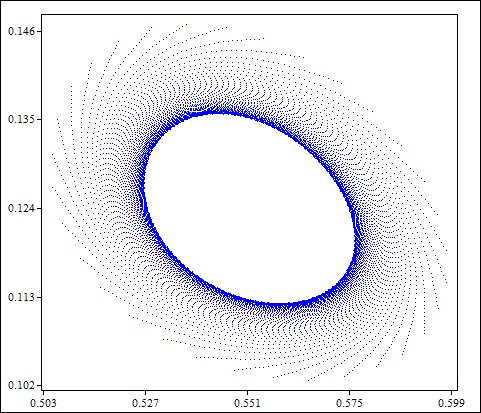}}\\
    \caption{Phase Portrait for the System \eqref{R3}}
    \label{PP}
\end{figure}
\end{example}

\begin{example}
For $\alpha = 0.5, h = 1, r = 2, d = 2.5$, $d=2.91667$, $x_0=0.8$ and $y_0=0.09$ then the system \eqref{R3} has a positive fixed point $(0.857143, 0.0979592)$. The Jacobian of the system \eqref{R3} at this fixed point has the eigenvalues $\lambda_1=-1$ and $\lambda_2=-0.428571$. Hence, system \eqref{R3} observes a period doubling bifurcation.
\begin{figure}
    \centering
    \subfigure[$(b, x_n)$]{\includegraphics[width=0.45\textwidth]{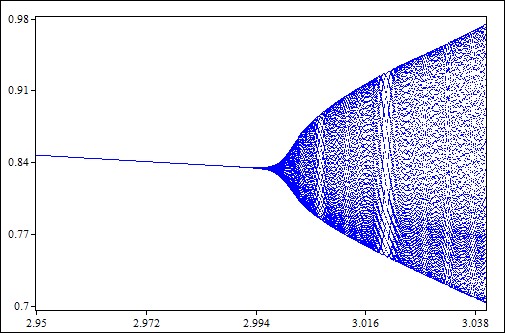}}
    \subfigure[$(b, y_n)$]{\includegraphics[width=0.45\textwidth]{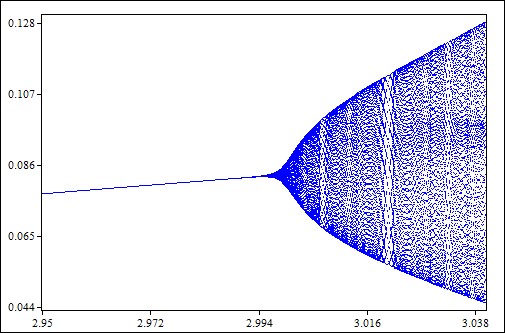}}\\
    \caption{Bifurcation Diagrams for the System \eqref{R3}}
    \label{BC}
\end{figure}
\end{example}

\section{Conclusion}
In \cite{XD}, the author uses the forward Euler scheme to discretize the system \eqref{R1} and then investigates the complex dynamics of the discrete-time predator-prey system \eqref{R1} in the closed first quadrant $R_{+}^2$. He observed that the system undergoes the flip bifurcation and Hopf bifurcation in the interior of the $R_{+}^2$ using the centre manifold theorem and bifurcation theory. He also presented the numerical simulation not only to illustrate the theoretical analysis but also to exhibit the complex dynamic behaviours, such as the period-5,6,9,10,14,18,20,25 orbits, a cascade of period-doubling bifurcation in period-2,4,8, quasi-periodic orbits and the chaotic sets. However, these results reveal much better results than those of the continuous system. \\
In our study, we use the piecewise constant argument technique, conformal fractional case, to discretize the system \eqref{R1} and then investigate its complex dynamics. We also presented numerical examples to illustrate our theoretical work. We observe that it gives more accurate results than those in \cite{XD}. Graphical illustrations are also presented to support numerical work.\\
\\
\textbf{Acknowledgement}
This research was funded in whole by the National Science Centre, Poland, grant number 2023/51/B/ST8/01062. For the purpose of Open Access, the authors have applied a CC-BY public copyright license to any Author Accepted Manuscript (AAM) version arising from this submission.\\
\\
\textbf{Declaration of Conflicts of Interest}
The authors affirm they have no conflicts of interest to disclose concerning the current study.\\
\\
\textbf{Author Contribution Statement} 
All authors contributed equally. Also, all authors reviewed the manuscript.\\
\\
\textbf{Availability of Data and Materials} 
The datasets analyzed during the current study are available from the corresponding author upon reasonable request.\\

\end{document}